\documentclass[a4paper,10pt]{amsart}

\usepackage{amssymb} \usepackage[OT4]{fontenc}
\usepackage{framed} \usepackage{graphicx}
\usepackage{pdfpages}
\usepackage{tikz}

\theoremstyle{definition}
\newtheorem{definition}{Definition}
\theoremstyle{plain} \newtheorem{theorem}[definition]{Theorem}

 \theoremstyle{remark} \newtheorem{remark}[definition]{Remark}

\newtheorem{problem}[definition]{Problem}

\def\P{\mathbb{P}}

\def\K{\mathbb{K}}
\def\R{\mathbb{R}}

\newcommand\call{{\mathcal L}}

\renewcommand\gamma{d}

\newcommand\rd[1]{\lfloor#1\rfloor}

\begin{document}

\title{A counterexample to the containment $I^{(3)}\subset I^2$ over the reals}

\author{Adam Czapli\'nski}
\author{Agata G\l ówka}
\author{Grzegorz Malara}
\author{Magdalena Lampa--Baczy\'nska}
\author{Patrycja \L uszcz--\'Swidecka}
\author{Piotr Pokora}
\author{Justyna Szpond}




\begin{abstract}
   The purpose of this note is to give defined over the real numbers counterexamples to a question
   relevant in the commutative algebra, concerning a containment relation between algebraic
   and symbolic powers of homogeneous ideal.
\end{abstract}

\maketitle
\section{The main result}
   In algebraic geometry and in commutative algebra there has been recently
   a lot of interest in comparing usual (algebraic) and symbolic powers of homogeneous
   ideals, see for example \cite{BocHar10}, \cite{HaHu13}, \cite{BCH14}.
   If $I\subset\K[x_0,x_1,\dots,x_n]$ is a homogeneous ideal, then its \emph{algebraic}
   $r$--th power $I^r$ is defined as the ideal generated by $r$--th powers $f^r$
   of all elements $f$ in $I$. This is a purely \emph{algebraic} concept.
   On the other hand, homogeneous ideals are defined geometrically as sets of
   all polynomials vanishing along a given set, a subvariety $V\subset\P^n(\K)$.
   For example, the ideal determined by simultaneous vanishing in points
   $P=(1:0:0)$, $Q=(0:1:0)$ and $R=(0:0:1)$ in the projective plane $\P^2(\K)$
   is generated by the monomials $yz, xz$ and $xy$. Studying the geometry of
   algebraic varieties one is often interested in polynomials vanishing to certain
   order $m$ along the subvariety. For example the monomial $xy$ vanishes to order
   $2$ in point $R$ and only to order $1$ in points $P$ and $Q$, whereas
   the monomial $xyz$ vanishes to order $2$ in all these points.
   All polynomials vanishing to order $m$ along the subvariety defined
   by an ideal $I$ form again an ideal. This ideal is denoted by $I^{(m)}$
   and is called the $m$--th \emph{symbolic} power of $I$. This is a \emph{geometric}
   concept. It follows from the definition that $I^m\subset I^{(m)}$
   holds for all homogeneous ideals. The reverse inclusion however might
   fail already in the simplest situations. For example this happens
   for the ideal $I$ of the three points $P, Q$ and $R$ defined above. 
   Indeed, the least
   degree of a polynomial in $I^2$ is $4$ (because generators of $I$ have degree $2$),
   whereas $xyz$ is a polynomial of degree $3$ contained in $I^{(2)}$.
   Hence there is surely no containment $I^{(2)}\subset I^2$ in this case.
   
   It came as a big surprise that nevertheless there
   is also a uniform containment relation in the reverse direction, when
   one takes into account the dimension of the ambient space (or equivalently
   the number of variables in the polynomial ring). More precisely,
   it has been discovered independently by Ein, Lazarsfeld and Smith \cite{ELS00}
   in characteristic $0$ and Hochster and Huneke \cite{HoHu02} in finite
   characteristic, that there is always the containment
   $$I^{(m)}\subset I^r,$$
   provided $m\geq n\cdot r$, where $n$ is the dimension of the ambient space.
   Whereas the lower bound $n\cdot r$ cannot be improved in general, it has
   been expected that it can be improved under certain assumptions on
   varieties defined by $I$. In particular, for points in the projective plane
   $\P^2(\K)$ Huneke asked if there is always the containment
   \begin{equation}\label{eq:cont32}
      I^{(3)}\subset I^2.
   \end{equation}
   It is the first instance of a more general statement predicting
   for ideals defined by vanishing along points the
   containment
   \begin{equation}\label{eq:cont gen}
      I^{(nr-(n-1))}\subset I^r
   \end{equation}
   for all $r\geq 1$.

   A first counterexample to the containment in \eqref{eq:cont32} was announced
   in 2013 by Dumnicki, Szemberg and Tutaj-Gasi\'nska in \cite{DST13}. In that
   counterexample the set of relevant points is taken as all the intersection
   points of $9$ lines whose union is defined by the polynomial equation
   $$(x^3-y^3)(y^3-z^3)(z^3-x^3)=0.$$
   These lines form the so called \emph{dual Hesse configuration}, see \cite{ArtDol09}
   for a beautiful account on that, interesting in its own right, subject.
   The $9$ lines in this configuration intersect by $3$ in $12$ points. Since
   there are no points where only two of these lines intersect, such a configuration
   cannot be realized over the reals as it would contradict the following
   well-know theorem
   of Sylvester \cite{Syl1893} and Gallai.
\begin{theorem}[Sylvester--Gallai]\label{thm:SG}
   Given a finite set of points in the real plane,
   either all points are collinear or there exists a pair of points
   not collinear with any other point in the set.
\end{theorem}
   The contradiction arises if one passes to the dual statement.
\begin{theorem}[dual Sylvester--Gallai Theorem]\label{thm:dSG}
   Given a finite number of lines in the real plane,
   either all these lines belong to the same pencil of lines
   passing through a fixed point, or there exists a point
   in the plane, where only two of these lines intersect.
\end{theorem}

   Shortly after the counterexample in \cite{DST13} was announced,
   Harbourne and Seceleanu \cite{HS13} constructed a series
   of counterexamples to the containment in \eqref{eq:cont gen}
   for various values of $n$ and $r$. Their counterexamples however are
   defined over finite fields.

   Since configurations of real lines are subject to stronger
   combinatorial constrains than sets of lines either in the complex
   plane or in planes defined over finite fields, it is not
   immediately clear that counterexamples to \eqref{eq:cont32}
   can be constructed in the real plane. This issue has motivated
   our research. The main result we want to announce is the following.
\begin{theorem}\label{thm:main}
   There exists a series of counterexamples to the containment in \eqref{eq:cont32}
   defined in the real plane.
\end{theorem}
   In fact our counterexamples are well known in the combinatorics. They were
   introduced by F\"uredi and Pal\'asti in \cite{FP84} in connection with
   the following Sylvester-Gallai problem motivated by Theorem \ref{thm:dSG}.
\begin{problem}
   Given a configuration of $s$ lines in the real (projective) plane,
   not all belonging to the same pencil of lines through a fixed point,
   what is the minimal number of points where only two lines meet?
\end{problem}
   It has been long conjectured, see for example \cite{BGS74},
   that there are at most $R_s=1+\rd{\frac{s(s-3)}{6}}$
   points where three lines meet, so that there are at least
   $\binom{s}{2}-3-3\rd{\frac{s(s-3)}{6}}$ points where the
   lines meet in pairs only. This conjecture, at least for large $s$,
   was recently proved by Green and Tao \cite{GreTao13}.
   A series of examples with number of triple points
   equal or close to $R_s$ was constructed in \cite{FP84}.

   In the recent preprint \cite{Tao13} Terence Tao pointed out
   that many combinatorial problems have been solved, or substantial
   progress has been obtained, when applying methods from algebra
   or algebraic geometry. In this note the directions are reversed,
   we present a solution to an algebraic problem based on a combinatorial
   construction.
\section{Proof of the main result}
   In order to prove Theorem \ref{thm:main}, we are interested in configurations
   of lines with a high number of points where three of them meet.
   The motivation for this approach
   is the following. Taking $V$ as the set of all triple points
   in the configuration and denoting by $I=I(V)$ the ideal of $V$,
   it is immediate from the description of symbolic powers given above,
   that the product
   $$F=L_1\cdot\ldots\cdot L_s$$
   of equations of the lines in the configuration is contained in $I^{(3)}$.

   If the number of triple points is high when related to the number
   of lines, then $F$ will be the only element of degree $s$ in $I^{(3)}$.
   On the other hand, the generators of $I$ should be of relatively high
   degree. In the most optimistic situation, the least degree $\gamma$
   of a generator of $I$ would satisfy $2\gamma>s$. In that case, it would
   be clear that $F$ cannot be contained in $I^2$. This scenario never
   happens in the reality, so that the actual argument requires somewhat
   more effort.

   We recall now how the examples of F\"uredi and Pal\'asti work.
   To this end it is convenient to identify
   the real plane $\R^2$ with the set of complex numbers $\mathbb{C}$
   in the usual way. Let $n$ be an even and positive integer.
   Let $\xi=\exp(2\pi i/n)$ be a primitive $n$--th root of unity
   and let $P_i=\xi^{i}$ for $i=0,\ldots, n-1$.
   Denote by $L_{m,k}$ the (real) line passing through the points $P_{m}$ and $P_{k}$ if $m\neq k$
   and the tangent line to the unit circle at the point $P_{m}$ if $m=k$.

   Then we construct
   the configuration of points and lines in two steps.
   First we define the configuration of lines
   $$\call_n = \{L_{i,\frac{n}{2}-2i} \}_{i=0}^{n-1},$$
   where the indices are understood ${\rm mod} \, n$.
   We will denote the linear form defining the line $L_{i,\frac{n}{2}-2i}$ by $L_i$ for $i=0,\ldots,n-1$.
   Let $F_n=\prod_{i=0}^{n-1}L_i$. Then $\deg(F_n)=n$.
   Let $Z_{n}$ be the set of all triple points in the configuration $\call_n$.

   Note that the lines $L_i, L_j, L_k$ are concurrent if and only if $n$ divides $i+j+k$, see \cite[Lemma]{FP84}.
   It follows from \cite[Property 4]{FP84} that there are exactly $1 + \lfloor n(n-3)/ 6 \rfloor $
   points in the set $Z_n$.

   By the way of an example we focus now on the case $n=12$. This is the
   minimal number of lines, which leads to a configuration of points
   giving a counterexample to the containment in \eqref{eq:cont32}.

   The configuration of lines arising in this case is illustrated
   on the following picture prepared with the aid of Geogebra.
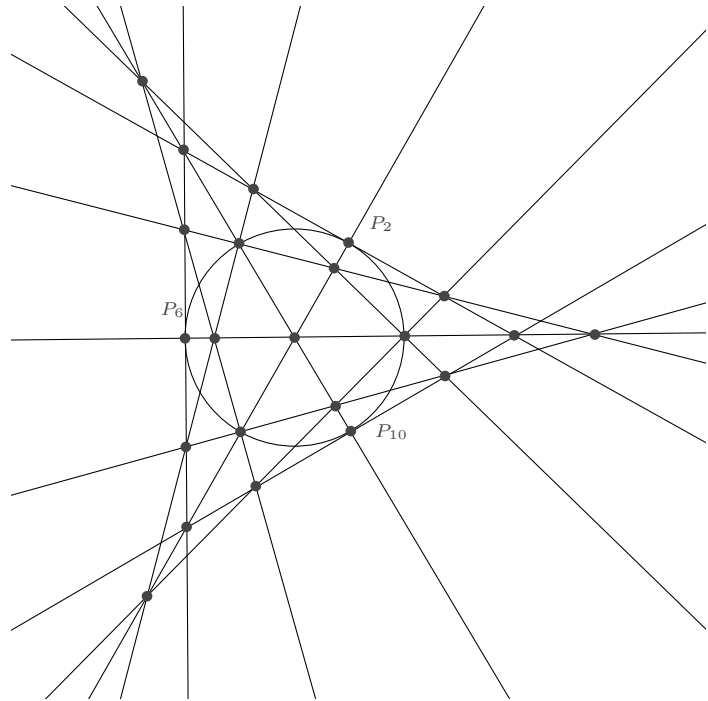
\begin{center}
\definecolor{uuuuuu}{rgb}{0.27,0.27,0.27}
\begin{tikzpicture}[line cap=round,line join=round,x=1.0cm,y=1.0cm]
\clip(-0.82,-3.96) rectangle (8.3,5.22);
\draw [domain=-4.36:18.32] plot(\x,{(-2.33-0.03*\x)/-2.89});
\draw [domain=-4.36:18.32] plot(\x,{(-5.21--0.51*\x)/-1.98});
\draw [domain=-4.36:18.32] plot(\x,{(--7.45-1.43*\x)/1.46});
\draw [domain=-4.36:18.32] plot(\x,{(--8.44-2.49*\x)/1.47});
\draw [domain=-4.36:18.32] plot(\x,{(--4.08-1.97*\x)/0.55});
\draw [domain=-4.36:18.32] plot(\x,{(-3.23--1.98*\x)/0.51});
\draw [domain=-4.36:18.32] plot(\x,{(-6.11--2.51*\x)/1.42});
\draw [domain=-4.36:18.32] plot(\x,{(-5.12--1.46*\x)/1.43});
\draw [domain=-4.36:18.32] plot(\x,{(--2.03-0.55*\x)/-1.97});
\draw [domain=-4.36:18.32] plot(\x,{(--5.19-0.71*\x)/1.26});
\draw [domain=-4.36:18.32] plot(\x,{(--3.18-0.73*\x)/-1.24});
\draw [domain=1:2] plot(\x,{(-2.12--1.44*\x)/-0.01});
\draw(2.9,0.83) circle (1.44cm);
\begin{scriptsize}
\fill [color=uuuuuu] (2.19,-0.42) circle (2.0pt);
\fill [color=uuuuuu] (4.35,0.85) circle (2.0pt);
\fill [color=uuuuuu] (2.17,2.08) circle (2.0pt);
\fill [color=uuuuuu] (0.9,4.23) circle (2.0pt);
\fill [color=uuuuuu] (1.44,3.32) circle (2.0pt);
\fill [color=uuuuuu] (1.45,2.26) circle (2.0pt);
\fill [color=uuuuuu] (2.36,2.8) circle (2.0pt);
\fill [color=uuuuuu] (1.85,0.82) circle (2.0pt);
\fill [color=uuuuuu] (2.9,0.83) circle (2.0pt);
\fill [color=uuuuuu] (1.48,-1.68) circle (2.0pt);
\fill [color=uuuuuu] (0.96,-2.6) circle (2.0pt);
\fill [color=uuuuuu] (2.39,-1.14) circle (2.0pt);
\fill [color=uuuuuu] (1.47,-0.62) circle (2.0pt);
\fill [color=uuuuuu] (3.44,-0.08) circle (2.0pt);
\fill [color=uuuuuu] (4.88,0.32) circle (2.0pt);
\fill [color=uuuuuu] (4.87,1.38) circle (2.0pt);
\fill [color=uuuuuu] (5.79,0.86) circle (2.0pt);
\fill [color=uuuuuu] (6.85,0.87) circle (2.0pt);
\fill [color=uuuuuu] (3.42,1.75) circle (2.0pt);
\fill [color=uuuuuu] (3.61,2.09) circle (2.0pt);
\draw[color=uuuuuu] (4.04,2.34) node {$P_2$};
\fill [color=uuuuuu] (3.64,-0.41) circle (2.0pt);
\draw[color=uuuuuu] (4.18,-0.42) node {$P_{10}$};
\fill [color=uuuuuu] (1.46,0.82) circle (2.0pt);
\draw[color=uuuuuu] (1.3,1.2) node {$P_6$};
\end{scriptsize}
\end{tikzpicture}
\end{center}
   The points in $Z_{12}$ are marked by dots.
   There are $19$ of them.
   There are $3$ lines, each containing $4$ configuration points, these lines are the tangents
   to the unit circle at points $P_2, P_6$ and $P_{10}$.
   On each of the remaining $9$ configuration lines there are exactly $5$ configuration points.

   Thus the configuration is in fact the union of two configurations.
   One depicted in the picture below using solid lines
   is a real analogue of the dual Hesse
   configuration mentioned before.
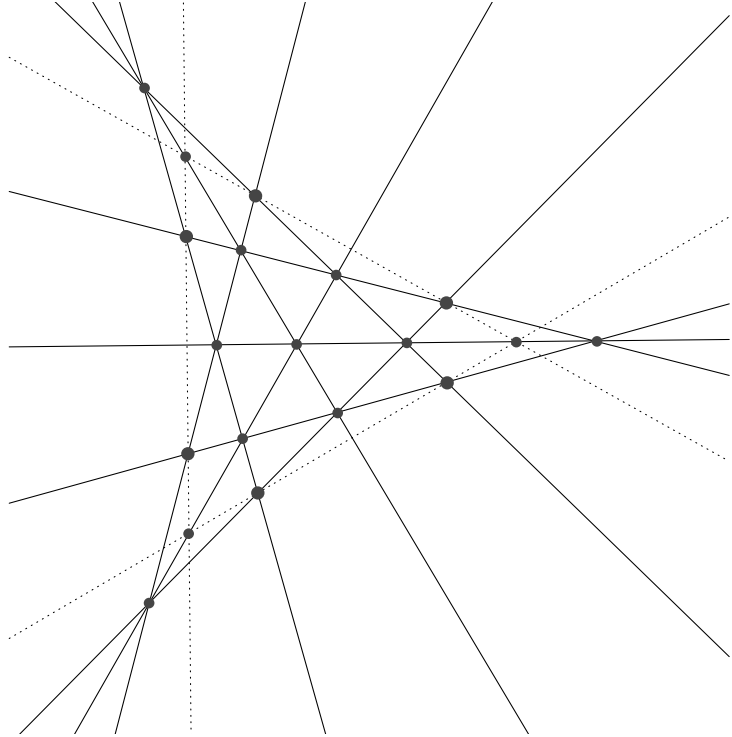
\begin{center}
\definecolor{uuuuuu}{rgb}{0.27,0.27,0.27}
\begin{tikzpicture}[line cap=round,line join=round,x=1.0cm,y=1.0cm]
\clip(-1.33,-4.33) rectangle (9.1,5.36);
\draw [domain=-0.88:8.59] plot(\x,{(-2.33-0.03*\x)/-2.89});
\draw [domain=-0.88:8.59] plot(\x,{(-5.21--0.51*\x)/-1.98});
\draw [domain=-0.88:8.59] plot(\x,{(--7.45-1.43*\x)/1.46});
\draw [domain=-0.88:8.59] plot(\x,{(--8.44-2.49*\x)/1.47});
\draw [domain=-0.88:8.59] plot(\x,{(--4.08-1.97*\x)/0.55});
\draw [domain=-0.88:8.59] plot(\x,{(-3.23--1.98*\x)/0.51});
\draw [domain=-0.88:8.59] plot(\x,{(-6.11--2.51*\x)/1.42});
\draw [domain=-0.88:8.59] plot(\x,{(-5.12--1.46*\x)/1.43});
\draw [domain=-0.88:8.59] plot(\x,{(--2.03-0.55*\x)/-1.97});
\draw [dotted,domain=1:2] plot(\x,{(--4.24-2.89*\x)/0.03});
\draw [dotted,domain=-0.88:8.59] plot(\x,{(--10.39-1.42*\x)/2.51});
\draw [dotted,domain=-0.88:8.59] plot(\x,{(--6.35-1.47*\x)/-2.49});
\begin{scriptsize}
\fill [color=uuuuuu] (2.19,-0.42) circle (2.0pt);
\fill [color=uuuuuu] (4.35,0.85) circle (2.0pt);
\fill [color=uuuuuu] (2.17,2.08) circle (2.0pt);
\fill [color=uuuuuu] (0.9,4.23) circle (2.0pt);
\fill [color=uuuuuu] (1.44,3.32) circle (2.0pt);
\fill [color=uuuuuu] (1.45,2.26) circle (2.5pt);
\fill [color=uuuuuu] (2.36,2.8) circle (2.5pt);
\fill [color=uuuuuu] (1.85,0.82) circle (2.0pt);
\fill [color=uuuuuu] (2.9,0.83) circle (2.0pt);
\fill [color=uuuuuu] (1.48,-1.68) circle (2.0pt);
\fill [color=uuuuuu] (0.96,-2.6) circle (2.0pt);
\fill [color=uuuuuu] (2.39,-1.14) circle (2.5pt);
\fill [color=uuuuuu] (1.47,-0.62) circle (2.5pt);
\fill [color=uuuuuu] (3.44,-0.08) circle (2.0pt);
\fill [color=uuuuuu] (4.88,0.32) circle (2.5pt);
\fill [color=uuuuuu] (4.87,1.38) circle (2.5pt);
\fill [color=uuuuuu] (5.79,0.86) circle (2.0pt);
\fill [color=uuuuuu] (6.85,0.87) circle (2.0pt);
\fill [color=uuuuuu] (3.42,1.75) circle (2.0pt);
\end{scriptsize}
\end{tikzpicture}
\end{center}
   In the real case there are only $10$ triple points. The remaining
   $6$ double points (distinguished in the picture by larger dots)
   determine in pairs, $3$ additional lines (dotted lines in the picture), which
   intersect in pairs on the lines of the solid configuration. These $3$ lines
   form the second configuration, and altogether there are $12$ lines and $19$ triple
   points.

   In order to determine generators of the ideal $J=I(Z_{12})$, we can intersect
   the ideals of single points in the set $Z_{12}$, i.e.
   $$J=I(Q_1)\cap\ldots\cap I(Q_{19}).$$
   Note that the ideals $I(Q_j)$ are generated by any two lines
   passing through the point $Q_j$, so they can be easily written
   explicitly down. In fact, these calculations can be done by hand,
   but this is a tedious procedure. Instead, we prefer to revoke
   a Singular \cite{DGP11} script in this place.

   In fact, our computational strategy in the script attached at the end of this note
   is slightly different (and more efficient). We compute first
   the points $P_0,\ldots,P_{11}$ on the unit circle, the vertices
   of the regular $12$--gon with $P_0=(1,0)$ (in the script we use projective
   coordinates, so that $P_0=(1:0:1)$). Then we compute the equations of lines
   $L_0,\ldots,L_{11}$. In the next step, we don't compute the triple
   points in the configuration $\call_{12}$ one by one. Instead, we
   compute the ideal $J$ of all of them simultaneously taking into
   account that they are the worst singularities of the configuration
   (in algebraic geometry a triple point is considered to be a worse singularity than
   a double point). Finally, we check the containment of $F_{12}$ (which is
   called $F$ in the script) in $J^2$. Again, this step could be
   carried out by hand in the theory, but this is a lengthy and dull
   computation, which can be safely relegated to a machine.

   Checking any given configuration with an even number $n\geq 12$
   goes along the same lines.

   This finishes the proof of Theorem \ref{thm:main}.
\endproof
\begin{remark}
   There is a similar construction for an odd number of lines $n$
   in \cite{FP84}. We have checked also these examples and found
   out that they give counterexamples to the containment
   in \eqref{eq:cont32} for all $n\geq 13$. Thus $12$ is the minimal
   number of real lines such that triple points of their configuration
   provide a counterexample to \eqref{eq:cont32}.
\end{remark}

\begin{remark}
   All counterexamples to the containment \eqref{eq:cont gen} found up to now
   are based on configurations of points coming from some extremal, from the
   combinatorial point of view, configurations of lines. It would be very
   interesting to understand if there are some deeper connections justifying
   these phenomena lurking behind the scenes.
\end{remark}

\section{A Singular script}
\begin{verbatim}
LIB "elim.lib";
ring R=(0,a),(x,y,z),dp; option(redSB);
minpoly=a2-3;
proc gline(ideal I1, ideal I2) {
  ideal I=intersect(I1,I2);
  I=std(I);
  return(I[1]);
}
ideal P0=x-z,y; ideal P1=x-a/2*z,y-1/2*z; ideal P2=x-1/2*z,y-a/2*z;
ideal P3=x,y-z; ideal P4=x+1/2*z,y-a/2*z; ideal P5=x+a/2*z,y-1/2*z;
ideal P6=x+z,y; ideal P7=x+a/2*z,y+1/2*z; ideal P8=x+1/2z,y+a/2*z;
ideal P9=x,y+z; ideal P10=x-1/2*z,y+a/2*z; ideal P11=x-a/2*z,y+1/2*z;

poly L0=gline(P0,P6); poly L1=gline(P1,P4); poly L2=1/2*x+a/2*y-z;
poly L3=gline(P3,P0); poly L4=gline(P4,P10); poly L5=gline(P5,P8);
poly L6=-x-z; poly L7=gline(P7,P4); poly L8=gline(P8,P2);
poly L9=gline(P9,P0); poly L10=1/2*x-a/2*y-z; poly L11=gline(P11,P8);

poly F=L0*L1*L2*L3*L4*L5*L6*L7*L8*L9*L10*L11;
ideal I=F,diff(F,x),diff(F,y),diff(F,z);
ideal J=I,diff(I,x),diff(I,y),diff(I,z);
ideal M=x,y,z; J=std(J);
J=sat(J,M)[1]; J=std(J);
reduce(F,std(J^2));
quit;
\end{verbatim}

\paragraph*{\emph{Acknowledgement.}}
   This note resulted from the second workshop on linear series held in Lanckorona
   in September 2013. We thank Pedagogical University of Cracow for financial
   support. Last but not least we thank Marcin Dumnicki, Tomasz Szemberg
   and Halszka Tutaj-Gasi\'nska for valuable comments.

\bigskip \small

\bigskip
   Agata G\l ówka--Habura, Grzegorz Malara, Magdalena Lampa--Baczy\'nska, Patrycja \L uszcz--\'Swidecka, Piotr Pokora, Justyna Szpond,
   Pedagogical University of Cracow, Institute of Mathematics,
   Podchor\c a\.zych 2,
   PL-30-084 Krak\'ow, Poland

\bigskip
    Adam Czapli\'nski,
    Johannes Gutenberg-Universit\"at Mainz, Institute of Mathematics,
    Staudingerweg 9,
    55099 Mainz

\end{document}